\documentclass[11pt,reqno]{amsart}

\usepackage{geometry}
\usepackage{url}

\usepackage{amsthm,amsmath,amssymb,bm,color,wrapfig}
\usepackage{ascmac}
\usepackage{graphicx}
\usepackage{multirow}
\usepackage[export]{adjustbox}
\usepackage[hidelinks]{hyperref}

\newcommand{\E}[0]{\mathbb{E}}

\newcommand{\R}[0]{\mathbb{R}}

\renewcommand{\Pr}[0]{\mathbb{P}}

\newcommand{\calR}{\mathcal{R}}
\newcommand{\calA}{\mathcal{A}}
\newcommand{\calB}{\mathcal{B}}
\newcommand{\calC}{\mathcal{C}}
\newcommand{\calF}{\mathcal{F}}
\renewcommand{\hat}{\widehat}
\renewcommand{\le}{\leqslant}
\renewcommand{\ge}{\geqslant}

\newcommand{\code}[1]{\texttt{#1}}

\DeclareMathOperator{\Var}{Var}

\DeclareMathOperator*{\argmin}{arg\,min}

\definecolor{fgcolor}{rgb}{0.345, 0.345, 0.345}

\usepackage{framed}
\makeatletter
 {\par\unskip\endMakeFramed%
 \at@end@of@kframe}
\makeatother

\definecolor{shadecolor}{rgb}{.97, .97, .97}
\definecolor{messagecolor}{rgb}{0, 0, 0}
\definecolor{warningcolor}{rgb}{1, 0, 1}
\definecolor{errorcolor}{rgb}{1, 0, 0}

\usepackage{listings}
\usepackage{alltt}

\theoremstyle{plain}
\newtheorem{theorem}{Theorem}[section]
\newtheorem{lemma}{Lemma}[section]
\newtheorem{proposition}{Proposition}[section]
\newtheorem{corollary}{Corollary}

\theoremstyle{definition}

\newtheorem{example}{Example}[section]
\newtheorem{remark}{Remark}

\begin{document}


\title[High-dimensional Data Bootstrap]{High-dimensional Data Bootstrap}

\thanks{V. Chernozhukov is partially supported by Amazon Core AI research grant.  
K. Kato is partially supported by NSF grants DMS-1952306 and DMS-2014636. 
Y. Koike is partially supported by JST CREST}
 \date{First version: April 11, 2022. This version: \today}

\author[V. Chernozhukov]{Victor Chernozhukov}

\address[V. Chernozhukov]{
Department of Economics and Center for Statistics and Data Science, MIT.
}
\email{vchern@mit.edu}

\author[D. Chetverikov]{Denis Chetverikov}

\address[D. Chetverikov]{
Department of Economics, UCLA.}
\email{chetverikov@econ.ucla.edu}

\author[K. Kato]{Kengo Kato}

\address[K. Kato]{
Department of Statistics and Data Science, Cornell University.
}
\email{kk976@cornell.edu}

\author[Y. Koike]{Yuta Koike}

\address[Y. Koike]{
Mathematics and Informatics Center and Graduate School of Mathematical Sciences, The University of Tokyo.}
\email{kyuta@ms.u-tokyo.ac.jp}

\begin{abstract}
This article reviews recent progress in high-dimensional bootstrap. We first review high-dimensional central limit theorems for distributions of sample mean vectors over the rectangles, bootstrap consistency results in high dimensions,  and key techniques used to establish those results. We then review selected applications of high-dimensional bootstrap: construction of simultaneous confidence sets
for high-dimensional vector parameters, multiple hypothesis testing via stepdown, post-selection inference, intersection bounds for partially identified parameters, and inference on best policies in policy evaluation. Finally, we also comment on a couple of future research directions. 
\end{abstract}


\maketitle


\section{Introduction}
\subsection{Overview}

The bootstrap is a generic method to estimate the sampling distribution of a statistic, typically by resampling one's own data. Since the seminal work of \cite{Efron1979}, there has been a substantial amount of research that explores the theoretical properties of the bootstrap. In classical settings where the data dimension is fixed, the bootstrap often yields more accurate confidence intervals or tests than those based on first-order asymptotic approximations. Also, the bootstrap provides a practical method of inference for models where the analytical estimation of the asymptotic distribution of an estimator or a test statistic is difficult, such as in quantile regression or semiparametric models.

Recently, there has been growing interest in extending the scope of the bootstrap to high-dimensional (a.k.a. ``$p \gg n$") settings. Due to the advancement of science and technology, data sets in which the number of features (e.g., genes) well exceeds the sample size (e.g., the number of patients) have become common in many application domains.
Analysis of such high-dimensional data has been a major focus in statistics in the last two or three decades; see, e.g., \cite{BuhlmannvandeGeer2011, Giraud2014, HastieTibshiraniWainwright2015, Wainwright2018} as textbook references on high-dimensional statistics. Yet, developing high-dimensional inference methods with provable accuracy under weak conditions is challenging, since classical statistical theory (which presumes fixed data dimensions) does not carry over to high dimensions, at least directly.

\cite{CCK2013AoS} contributed to the above goal and established the consistency of the Gaussian multiplier (or wild) and empirical bootstraps for the maximum of the sum of independent high-dimensional random vectors. Notably, \cite{CCK2013AoS} allow the data dimension $p$ to be \textit{much larger than} the sample size $n$. 
\cite{CCK2017AoP} extended the results of \cite{CCK2013AoS} and proved that, for $S_{n} = \sum_{i=1}^{n}X_{i}/\sqrt{n}$ the scaled sum of independent centered random vectors in $\R^{p}$ with $p = p_{n} \to \infty$, it holds that
\begin{equation}
\sup_{R \in \calR}\Big| \Pr\big(S_{n}^B \in R \mid X_1,\dots,X_n\big) - \Pr(S_n\in R)\Big| = O_P \left ( \frac{\log^7 (pn)}{n}\right )^{1/6}
\label{eq: CCK bound}
\end{equation}
under moment conditions, where $\calR$ is the class of rectangles in $\R^p$ and $S_n^B$ is either the Gaussian multiplier or empirical bootstrap statistic. Importantly, the error bound (\ref{eq: CCK bound}) only requires $p$ to be $\log p = o(n^{1/7})$, which allows $p \gg n$, for the bootstrap to be consistent.
 Subsequent works have explored refinements and extensions of these results, including \cite{deng2020beyond,lopes2020bootstrapping,kuchibhotla2021high,koike2021notes,chernozhukov2019improved,fang2021high,lopes2020central,kuchibhotla2020high,chernozhukov2020nearly}.

These extensive theoretical developments stimulated many new applications of bootstrap to high-dimensional or nonparametric statistics, including simultaneous inference for high-dimensional (or infinite-dimensional) parameters \cite{CCK2013AoS,CCK2014AoS,CCK2014AoSdensity,wasserman2014,Belloni2015,chen2015,chen2016,ZhangCheng2017JASA,Dezeure2017,chang2017simulation, ning2017general,belloni2018uniformly,rinaldo2019,
kuchibhotla2020valid}, testing for shape restrictions \cite{chetverikov2019testing}, detection of spurious correlations \cite{fan2018}, comparison of large covariance matrices \cite{chang2017}, inference for partially identified models  \cite{chetverikov2018adaptive,CCK2018RES},  goodness-of-fit testing \cite{jankova2020goodness}, error estimation of matrix randomized algorithms \cite{Lopes2019JMLR}, testing the mean function for functional data \cite{lopes2020bootstrapping}, and many more. 

This article aims to provide a brief overview of the current literature on high-dimensional bootstrap.   In Section \ref{sec: challenge}, we review classical asymptotics for the empirical bootstrap and discuss challenges in inference for high-dimensional data. In Section \ref{sec: main}, we review basic high-dimensional central limit theorems (CLTs) and bootstrap consistency results in high dimensions. In Section \ref{sec: app}, we discuss selected applications of high-dimensional bootstrap. 
Section \ref{sec: discussion} leaves some concluding remarks. 

Finally, while this review aims to convey the main ideas and techniques in high-dimensional bootstrap, the literature is now quite broad and has seen rapid expansion, and thus this review is by no means exhaustive. For instance, due to space limitation, we omit several topics such as  the empirical process extension \cite{CCK2014AoS,chernozhukov2016empirical}, high-dimensional CLTs and bootstrap for dependent data \cite{ZhangWu2017, ZhangCheng2017, CCK2018RES, chang2021central,chiang2021inference,kurisu2021gaussian} and $U$-statistics \cite{Chen2018,chen2019randomized,chen2020jackknife,song2019approximating,song2020stratified}.

\subsection{Notation}
We use $N(\mu,\Sigma)$ to denote the normal distribution with mean vector $\mu$ and covariance matrix $\Sigma$. Abusing the notation, we also use $N(\mu,\Sigma)$ to denote a random vector following $N(\mu,\Sigma)$. Let $\stackrel{d}{\to}$ and $\stackrel{P}{\to}$ denote convergence in distribution and convergence in probability, respectively. Let $\| \cdot \|$ and $\| \cdot \|_{\infty}$ denote the Euclidean and  $\ell^\infty$-norms for vectors, respectively,  i.e., for $x = (x_1,\dots,x_p)^T$, $\| x \| = \sqrt{\sum_{j=1}^p x_j^2}$ and $\| x \|_{\infty} = \max_{1 \le j \le p} |x_j|$.
The notation $\delta_x$ denotes the Dirac delta at point $x$. For two numbers $a,b \in \R$, we use the notation $a \vee b = \max \{ a,b \}$. 
We use $\lesssim$ to denote inequalities up to numerical constants. 

\section{Challenges in Inference for High-dimensional Data}
\label{sec: challenge}

\subsection{Classical Asymptotics}

We first review how the bootstrap works in the classical setting where the data dimension is fixed while the sample size $n$ tends to infinity. Consider first inference on the mean parameter $\mu$ from independent and identically distributed (i.i.d.) univariate ($p=1$) random variables $X_1,\dots,X_n$. As long as the common distribution has finite second moment, the CLT yields that, for the sample mean $\overline{X}_n = n^{-1}\sum_{i=1}^n X_i$,
\[
\sqrt{n}(\overline{X}_n - \mu) \stackrel{d}{\to} N(0,\sigma^2),
\]
where $\sigma^2 > 0$ is the population variance. We shall apply the empirical bootstrap to estimate the limit distribution $N(0,\sigma^2)$ and hence the sampling distribution of $\sqrt{n}(\overline{X}_n-\mu)$. To this end, conditionally on the data $X_1,\dots,X_n$,  generate an independent sample $X_1^B,\dots,X_n^B$ from the empirical distribution $P_n = n^{-1}\sum_{i=1}^n \delta_{X_i}$ of $X_1,\dots,X_n$. Since,  conditionally on the data, $P_n$ has mean $\overline{X}_n$ and variance $\hat{\sigma}^2 = n^{-1}\sum_{i=1}^n (X_i-\overline{X}_n)^2$, one can expect that, conditionally on the data,
\[
\sqrt{n}(\overline{X}_n^B - \overline{X}_n) \stackrel{d}{\approx} N(0,\hat{\sigma}^2)  \stackrel{d}{\approx} N(0,\sigma^2),
\]
where $\overline{X}_n^B = n^{-1}\sum_{i=1}^n X_i^B$.
This suggests that the conditional distribution of $\sqrt{n}(\overline{X}_n^B - \overline{X}_n)$ given the data approximates the $N(0,\sigma^2)$ distribution. Indeed, since $N(0,\sigma^2)$ has a continuous distribution function, it is not difficult to show that 
\[
\sup_{t \in \R} \Big | \Pr^B \big(\sqrt{n}(\overline{X}_n^B - \overline{X}_n) \le t\big) - \Pr \big(N(0,\sigma^2) \le t\big) \Big | \stackrel{P}{\to} 0. 
\]
Here and in what follows, $\Pr^B$ denotes the conditional probability given the data. This result implies that the conditional distribution of $\sqrt{n}(\overline{X}_n^B - \overline{X}_n)$ given the data consistently estimates the sampling distribution of $\sqrt{n}(\overline{X}_n-\mu)$.

Suppose, more generally, that we want to make inference on a scalar parameter $\theta^{\star}$, for which a suitable  estimator $\hat{\theta}_n$ (based on i.i.d. data $X_1,\dots,X_n$) is available. Suppose also that $\hat{\theta}_n$ admits an asymptotically linear expansion of the form
\begin{equation}
\label{eq: expansion}
\hat{\theta}_n - \theta^{\star} = \frac{1}{n}\sum_{i=1}^n \psi (X_i) + o_{P}(n^{-1/2}), 
\end{equation}
where $\psi$ is an influence function such that $\psi (X_i)$ has mean zero. This expansion implies that
$
\sqrt{n}(\hat{\theta}_n - \theta^{\star}) \stackrel{d}{\to} N(0,\sigma^2_\psi),
$
where $\sigma^2_{\psi}$ is the variance of $\psi (X_1)$ (assume  $\sigma^2_{\psi} > 0$). To estimate the $N(0,\sigma^2_\psi)$ distribution or the sampling distribution of $\sqrt{n}(\hat{\theta}_n - \theta^{\star})$, we may apply bootstrap to $\hat{\theta}_n$ by replacing the data $X_1,\dots,X_n$ with the bootstrap sample $X_1^B,\dots,X_n^B$, which gives rise to the bootstrap version of the estimator $\hat{\theta}_n^B$. Then, under regularity conditions, we can expect that the bootstrap estimator $\hat{\theta}_n^B$ admits a similar asymptotic linear expansion as
\begin{equation}
\label{eq: bootstrap expansion}
\hat{\theta}_n^B - \theta^{\star} = \frac{1}{n}\sum_{i=1}^n \psi (X_i^B) + o_{P}(n^{-1/2}).
\end{equation}
Subtracting (\ref{eq: expansion}) from (\ref{eq: bootstrap expansion}), we have
\[
\sqrt{n}(\hat{\theta}_n^B - \hat{\theta}_n) = \frac{1}{\sqrt{n}} \sum_{i=1}^n (\psi (X_i^B) - \overline{\psi}_n) + o_{P}(1),
\]
where $\overline{\psi}_n = n^{-1}\sum_{i=1}^n \psi (X_i)$. Thus, as in the sample mean case, we obtain the bootstrap consistency
\[
\sup_{t \in \R} \Big | \Pr^B \big(\sqrt{n}(\hat{\theta}_n^B - \hat{\theta}_n) \le t\big) - \Pr \big(N(0,\sigma_{\psi}^2) \le t\big) \Big | \stackrel{P}{\to} 0. 
\]
See, e.g., \cite{wellner1996bootstrapping}.

The bootstrap consistency result yields methods to construct (asymptotically valid) confidence intervals for $\theta^{\star}$. For example, let $\hat{q}_\alpha$ denote the conditional $\alpha$-quantile of $\hat{\theta}_n^B$ given the data, i.e., $\hat{q}_\alpha = \inf\{ t : \Pr^B\big( \hat{\theta}_n^B \le t\big) \ge \alpha \}$.
Then, from the bootstrap consistency result, it is not difficult to verify that the following percentile confidence interval, $[ 2\hat{\theta}_n - \hat{q}_{1-\alpha/2}, 2\hat{\theta}_n - \hat{q}_{\alpha/2}  ]$,
contains $\theta^{\star}$ with probability $1-\alpha+o(1)$ (cf. Lemma 23.3 in \cite{van2000asymptotic}). 
We refer the reader to \cite{Hall1992,VW1996,horowitz2001bootstrap} for several other classical asymptotic results on the bootstrap, including higher-order refinements and empirical process extension.

\subsection{Challenges in High Dimensions}

The discussion in the preceding section immediately extends to the case where $\overline{X}_n$ or $\hat{\theta}_n$ are multivariate, as long as the dimension is \textit{fixed}. However, in modern statistical applications, we are often interested in situations where the data dimension $p$ is comparable to or even much larger than the sample size $n$. Mathematically, such scenarios can be captured by allowing $p$ to depend on $n$ and considering the case where $p_n \to \infty$ as $n \to \infty$.

The main issue in high-dimensional inference problems is the lack of explicit limit distributions. To fix ideas, let $X_1,\dots,X_n$ be i.i.d. $p$-dimensional random vector with mean $\mu$ and covariance matrix $\Sigma$. Even when $p = p_n \to \infty$, one can still expect that $S_n = \sqrt{n}(\overline{X}_n - \mu)$ be approximated by $N(0,\Sigma)$, but since the dimension of $N(0,\Sigma)$ is $p = p_n \to \infty$, there will be no obvious ``final'' limit distribution for $S_n$. Of course, we can think of $S_n$ as random elements in $\R^\mathbb{N}$ (the countable product of $\R$), but weak convergence in $\R^{\mathbb{N}}$ is equivalent to finite dimensional convergence, which is too weak a result to develop genuinely high-dimensional inference methods.

Also, often, we are interested in approximating the sampling distribution of a certain functional of $S_n$, such as the $\ell^\infty$-norm, i.e., $\| S_n \|_{\infty} = \max_{1 \le j \le p} | S_{n,j}|$. In some situations, it could be possible to derive (after a proper normalization) limit distributions to such one-dimensional functionals; however, required regularity conditions tend to be restrictive. For instance, the $\ell^\infty$-functional $\| S_n\|_{\infty}$ may converge to the type I extreme value distribution after a proper normalization, but the derivation of such results typically requires the coordinates of the random vectors to be weakly dependent; see, e.g., Section 2.2 in \cite{chang2017} for relevant discussion.
Further, in practical situations, we often deal with high-dimensional statistics that are approximately linear (like (\ref{eq: expansion})), but how fast the linearization error should vanish is a nontrivial question because of the lack of explicit limit distributions.

Another challenge is that a large number of features induce a potentially complex dependence structure among the features. Failing to take into account the dependence structure among the features would result in (sometimes severely) conservative inference methods, especially in high dimensions. For instance, consider testing the global hypothesis $H_0: \mu_1=\cdots=\mu_p$ using $\| S_n \|_{\infty}$. Assume each coordinate of $S_n$ has unit variance (for simplicity) and apply a \v{S}idak-type correction, which yields a critical value $t_\alpha$ given by $\Pr (|Z| > t_\alpha) = 1-(1-\alpha)^{1/p}$ for $Z \sim N(0,1)$. This is the critical value computed as if the coordinates were independent  and each coordinate of $S_n$ were distributed as $N(0,1)$ \cite{fan2007many}. However, in the extreme situation where the coordinates are perfectly correlated, then the actual rejection probability under the null is, using the classical Berry-Esseen theorem, $\Pr(\| S \|_{\infty} > t_\alpha) = \Pr(|Z| > t_{\alpha}) + O(n^{-1/2}) =  1-(1-\alpha)^{1/p} + O(n^{-1/2}) = -p^{-1}\log(1-\alpha) + O(p^{-2} + n^{-1/2})$, which is very close to $0$ when $p$ is large, thus yielding conservative tests. Of course, the perfect correlation is an extreme case, but the above discussion indicates that if $S_n$ contains a large number of strongly correlated components, then Bonferoni or \v{S}idak-type corrections would yield overly conservative tests.

We will review below several techniques (high-dimensional CLTs, Gaussian anticoncentration inequality, Gaussian comparison) that help overcome the lack of explicit limit distributions and justify bootstrap methods for high-dimensional data. Also, we will demonstrate that the bootstrap is able to automatically take into account the dependence among the coordinates and yields asymptotically exact inference methods even in high dimensions. It should be noted that the bootstrap allows for arbitrary correlations between the coordinates, thereby accommodating a broader class of high-dimensional inference tasks.

\section{Bootstrap Methods in High Dimensions}
\label{sec: main}

A key step toward establishing the validity of bootstrap methods in high dimensions is a high-dimensional CLT. The high-dimensional CLT amounts to bounding a certain discrepancy measure between the sampling distribution of the scaled sample mean and the corresponding Gaussian distribution in such a way that the bound permits the dimension to increase with $n$. One natural approach in this direction is to compare the distributions over a (sufficiently large) collection of subsets $\calA$ in $\R^p$ and try to derive as explicit bounds as possible on $\sup_{A \in \calA} | \Pr (S_n \in A) - \Pr (N(0,\Sigma) \in A)|$, where $S_n$ is the scaled sample mean of mean-zero independent random vectors in $\R^p$ and $N(0,\Sigma)$ is the corresponding Gaussian distribution.
As \cite{CCK2013AoS,CCK2017AoP} observed, choosing $\calA$ to be the collection of rectangles allows us to establish error bounds that depend on the dimension $p$ only through $\log p$, thereby permitting $p \gg n$. Such high-dimensional CLTs and accompanying bootstrap consistency results, together with techniques developed therein, paved the way for constructing inference methods with theoretical guarantees for many high-dimensional inference problems under relatively mild regularity conditions.
This section reviews high-dimensional CLTs over the rectangles and  bootstrap consistency results in high dimensions, following the recent work of \cite{chernozhukov2019improved} that improves earlier results of \cite{CCK2013AoS,CCK2017AoP}. We also briefly review relevant developments in the literature.
Throughout this and the following sections, we always assume $n\geq2$ and $p\geq2$.

\subsection{High-Dimensional CLTs}

Let $X_1,\dots,X_n$ be independent (but not necessarily identically distributed) random vectors with dimension $p$. Assume that $X_i$'s have mean zero (otherwise, work with $X_i - \E[X_i]$ instead of $X_i$). Consider the scaled sample mean 
\[
S_n = \frac{1}{\sqrt{n}} \sum_{i=1}^n X_i. 
\]
Let $\overline{\sigma}, \underline{\sigma}$ be given positive constants such that $\underline{\sigma} \le \overline{\sigma}$, and let $B_n \ge 1$ be a sequence of constants that may diverge as $n \to \infty$. Let $\Sigma = \E[S_nS_n^T] = n^{-1}\sum_{i=1}^n \E[X_iX_i^T]$. Also,
let $\calR$ denote the collection of closed rectangles in $\R^p$,
\[
\calR = \Big \{ \prod_{j=1}^p [a_j,b_j] : -\infty \le a_j \le b_j \le \infty, \ j=1,\dots,p \Big \}.
\]
Finally, for notational convenience, let 
\[
\delta_{1,n} = \left(\frac{B_n^2\log^5(p n)}{n}\right)^{1/4} \quad \text{and} \quad 
\delta_{2,n}^{[q]} =\sqrt{ \frac{B_n^2(\log (pn))^{3-2/q}}{n^{1-2/q}}} \quad \text{for $q>2$}.
\]
We first present a high-dimensional CLT over the rectangles under a sub-exponential condition on the coordinates.

\begin{theorem}[High-dimensional CLT]\label{thm: HDCLT}
Suppose $\max_{1 \le i \le n; 1 \le j \le p} \E[e^{|X_{ij}|/B_n}] \le 2$, $\min_{1 \le j \le p} n^{-1}\sum_{i=1}^n \E[X_{ij}^2] \ge \underline{\sigma}^2$, and $\max_{1 \le j \le p} n^{-1}\sum_{i=1}^n \E[X_{ij}^4] \le B_n^2 \overline{\sigma}^2$. Then, 
\begin{equation}\label{eq: gaussian bound main}
\sup_{R \in \calR} \left|\Pr\left( S_n \in R \right) - \Pr (N(0,\Sigma) \in R) \right| \leq C\delta_{1,n},
\end{equation}
where $C$ is a constant that depends only on $\underline{\sigma}$ and $\overline{\sigma}$. 
\end{theorem}

The theorem follows from Lemma 4.3 in \cite{chernozhukov2019improved}. 
The assumption of the theorem is satisfied if, for example, $|X_{ij}| \le B_n$ almost surely for all $(i,j)$ and $\underline{\sigma}^2 \le n^{-1}\sum_{i=1}^n \E[X_{ij}^2] \le \overline{\sigma}^2$ for all $j$. Notably, the above theorem does not impose any restrictions on the correlation structure between the coordinates of the random vectors, so $\Sigma$ is permitted to be singular. This aspect is relevant in applications to high-dimensional data since, in the presence of many features, some of them are likely to have strong correlations, causing $\Sigma$ to be (or close to being) singular. For example, data from microarray or transcriptome experiments contain genes
that are divided into groups with varying sizes according
to functionalities, and genes from the same group tend to have
relatively strong (sometimes very strong) within-correlations (cf. \cite{chang2017}).

Theorem \ref{thm: HDCLT} shows that, under the assumption of the theorem, 
\begin{equation}
\sup_{R \in \calR} \left|\Pr\left( S_n \in R \right) - \Pr (N(0,\Sigma) \in R) \right| \to 0,
\label{eq: HDCLT}
\end{equation}
provided that $B_n^2 \log^5 (pn) = o(n)$, 
which allows $p$ to be much larger than $n$.

Theorem \ref{thm: HDCLT} requires that the coordinates of $X_i$ be sub-exponential. 
The following theorem, adapted from Theorem 2.5 in \cite{chernozhukov2019improved},  presents a version of the high-dimensional CLT under polynomial moment conditions.

\begin{theorem}[High-dimensional CLT under polynomial moment condition]
\label{thm:hdclt-poly}
Suppose that $\max_{1 \le i \le n}\E \left[ \| X_{i} \|_{\infty}^q \right]\leq B_n^q$ for some $q > 2$, $\min_{1 \le j \le p} n^{-1}\sum_{i=1}^n\E[X_{ij}^2] \ge \underline{\sigma}^2$, and $\max_{1 \le j \le p} n^{-1}\sum_{i=1}^n \E[X_{ij}^4] \le B_n^2 \overline{\sigma}^2$.
 Then,
\begin{align*}
\sup_{R \in \mathcal R} \left|\Pr\left( S_n \in R \right) - \Pr (N(0,\Sigma) \in R) \right|
\leq C (\delta_{1,n} \vee \delta_{2,n}^{[q]}),
\end{align*}
where $C$ is a constant that depends only on $q$, $\underline{\sigma}$ and $\overline{\sigma}$. 
\end{theorem}

Theorem \ref{thm:hdclt-poly} covers the following scenario relevant to regression applications: $X_i = \epsilon_i z_i$ where $\epsilon_i$ is a univariate ``error'' term while $z_i \in \R^p$ is a vector of fixed ``covariates''. In this case, $\E[\| X_i \|_{\infty}^q] \le \| z_i \|_{\infty}^q  \E[|\epsilon_i|^q]$, so if the covariates are uniformly bounded and the $q$-th moments of the error terms are bounded, then $B_n = O(1)$. Notably this only requires $\epsilon_i$ to have $q=2+ \delta$ bounded moments.

In practice, statistics of interest may not be exactly sample means but be only approximately linear. Since, in high dimensions, the approximating Gaussian distribution changes with $n$, the linearization error has to vanish at a sufficiently fast rate for such approximate sample means to satisfy the high-dimensional CLT. Fortunately, for the  rectangle case, the required condition on the linearization error is rather mild, as shown next. 

\begin{lemma}[High-dimensional CLT for approximate sample mean]
\label{lem: ASM}
Suppose that (\ref{eq: HDCLT}) holds for $S_n$, but $S_n$ is not directly available. 
Suppose instead that we have access to $\hat{S}_n$ that approximates $S_n$ such that $\hat{S}_n = S_n + \mathsf{R}_n$ with $\| \mathsf{R}_n \|_{\infty} = o_{P}(1/\sqrt{\log p})$. Assume $\min_{1 \le j \le p}\Sigma_{jj} \ge \underline{\sigma}^2$. Then the same conclusion (\ref{eq: HDCLT}) holds with $S_n$ replaced by $\hat{S}_n$. 
\end{lemma}

The lemma follows from the following Gaussian anticoncentration inequality for the rectangles, due (essentially) to \cite{nazarov2003maximal}; see also \cite{klivans2008learning}. A self-contained proof can be found in \cite{chernozhukov2017detailed}.
 \begin{proposition}[Nazarov's inequality]
 \label{lem: Nazarov}
Let $W = (W_1,\dots,W_p)^T$ be a (not necessarily centered) Gaussian random vector such that $\min_{1 \le j \le p} \Var (W_j) \ge \underline{\sigma}^2$ for some $\underline{\sigma} > 0$. Then, for $W^{\vee} = \max_{1 \le j \le p}W_j$, we have
\begin{equation}
\Pr \big(t < W^{\vee} \le t+\delta \big) \le  \frac{\delta}{\underline{\sigma} } (\sqrt{2\log p}+2), \  t \in \R, \delta > 0. 
\label{eq: nazarov}
\end{equation}
\end{proposition}

Nazarov's inequality states that the Gaussian maximum $W^{\vee}$ possesses a mass at most $O(\delta \sqrt{\log p})$ around the $\delta$-neighborhood of any point $t$. This is an opposite statement to the Gaussian concentration inequality, which asserts that the Gaussian maximum possesses a large mass around its expectation or median (cf. \cite{boucheron2013concentration}).
As such, Nazarov's inequality is an instance of an \textit{anticoncentration} inequality for the Gaussian distribution. The $\sqrt{\log p}$ dependence in Proposition \ref{lem: Nazarov} is sharp in the sense that the opposite inequality (up to multiplicative constants) holds when $W_1,\dots,W_p$ are i.i.d. as $N(0,1)$; see \cite{CCK2015PTRF}. Nazarov's inequality also plays a crucial role in establishing high-dimensional CLTs over the rectangles. 

To see how Lemma \ref{lem: ASM} follows from Nazarov's inequality, let us pick a particular rectangle $R = \{ x \in \R^p: \max_{1 \le j \le p} x_j \le t \}$ for some $t \in \R$ (the general case follows by translation and scaling; cf. the proof of Corollary 5.1 in \cite{CCK2017AoP}). Then, for $\hat{S}_n^{\vee} = \max_{1 \le j \le p}\hat{S}_{n,j}$ and $W \sim N(0,\Sigma)$, 
\[
\begin{split}
\Pr \big(\hat{S}_n \in R\big) &= \Pr \big(\hat{S}_n^{\vee} \le t\big) \le \Pr \big(S_n^{\vee} \le t + \epsilon\big) + \Pr (\| \mathsf{R}_n \|_{\infty} > \epsilon) \\
&=  \Pr \big(W^{\vee} \le t + \epsilon\big) + \Pr (\| \mathsf{R}_n \|_{\infty} > \epsilon) + o(1) \quad (\text{by (\ref{eq: HDCLT})}) \\
& = \Pr \big(W^{\vee} \le t\big) + O\big(\epsilon \sqrt{\log p}\big) + \Pr (\| \mathsf{R}_n \|_{\infty} > \epsilon) + o(1), \quad (\text{by (\ref{eq: nazarov})})
\end{split}
\]
where $o(1)$ is uniform in $t$. Likewise, we have
\[
\Pr \big(\hat{S}_n \in R\big) = \Pr \big(\hat{S}_n^{\vee} \le t\big) \ge \Pr \big(W^{\vee} \le t\big) - O\big(\epsilon \sqrt{\log p}\big) - \Pr (\| \mathsf{R}_n \|_{\infty} > \epsilon) - o(1).
\]
Since $\| \mathsf{R}_n \|_{\infty} = o_P(1/\sqrt{\log p})$, we can choose $\epsilon = \epsilon_n = o(1/\sqrt{\log p})$ in such a way that $\Pr (\| \mathsf{R}_n \|_{\infty} > \epsilon)  = o(1)$, which yields that $\sup_{t \in \R} |\Pr (S_n^{\vee} \le t) - \Pr (W^{\vee} \le t)| = o(1)$.

\begin{remark}[Gaussian anticoncentration]
The anticoncentration inequality in Proposition \ref{lem: Nazarov} is dimension-dependent in the sense that the bound depends explicitly on the dimension $p$. There are dimension-independent versions of anticoncentration inequalities for Gaussian maxima. For example, \cite{CCK2015PTRF} show that, if $W = (W_1,\dots,W_p)^T$ is centered Gaussian such that the coordinates have equal variance $\sigma^2 > 0$, then 
\[
\Pr\big (|W^{\vee} -t| \le \delta\big) \le 4 \delta \big(\E\big[W^{\vee}\big] + 1\big)/\sigma, \ \delta > 0.
\]
Since the latter bound does not explicitly depend on the dimension, the inequality readily extends to a centered Gaussian process $W = (W(t))_{t \in T}$ with homogeneous variance function such that $\sup_{t \in T} W(t) = \sup_{t \in T_0}W(t)$ almost surely for some countable subset $T_0$ of $T$. Indeed, the same inequality holds with $W^\vee$ replaced by $\sup_{t \in T}W(t)$ by passing $p$ to $\infty$. Such dimension-independent anticoncentration inequalities play an important role in nonparametric statistical inference. See \cite{CCK2014AoSdensity}. 
\end{remark}

\begin{remark}[Bibliographic notes]
Regarding the high-dimensional CLT over the rectangles, Theorem \ref{thm: HDCLT} is currently the sharpest bound available in the literature for sub-exponential coordinates under general correlation structures (i.e., without assuming nondegeneracy of $\Sigma$), both in terms of dependence on $n$ and $p$. \cite{CCK2013AoS} originally proved a version of the high-dimensional CLT for the maximum coordinate of $S_n$, namely, under the assumption of Theorem \ref{thm: HDCLT} (and one more technical moment condition):
\[
\sup_{A \in \calA} \Big|\Pr \big(S_n \in  A\big) - \Pr \big(N(0,\Sigma) \in A \big) \Big| \lesssim \left ( \frac{B_n^2 \log^b (pn) }{n} \right)^{1/a}
\]
for $\calA = \{ \{ x \in \R^p : \max_{1 \le j \le p} x_j \le t \} : t \in \R \} \subset \calR$ with $(a, b) =(8,7)$. Subsequently, \cite{CCK2017AoP} derived a version of Theorem \ref{thm: HDCLT} with
$(a, b) =(6,7)$ for $\calA = \calR$, which was further improved to $(a, b) =(6,5)$ for $\calA = \calR$ by \cite{koike2021notes}. Theorem \ref{thm: HDCLT} is taken from \cite{chernozhukov2019improved}, which improves these earlier results  to $(a, b) =(4,5)$.  See also \cite{kuchibhotla2021high} for non-uniform versions of high-dimensional CLTs. Finally, several recent works have established results with $a=2$ or near $2$, i.e.,
near $n^{-1/2}$ rates, under structural assumptions on $\Sigma$; see Section \ref{sec: parametric rate}. 

Theorems \ref{thm: HDCLT} and \ref{thm:hdclt-poly} are instances of multivariate Berry-Esseen bounds. 
There is a large and rich literature on this subject. Earlier references include \cite{asriev1986convergence,gotze1991rate,bentkus2003dependence,bentkus2005lyapunov,chatterjee2008multivariate,reinert2009multivariate}, among others.  \cite{CCK2013AoS} (Appendix L) provide a review of this literature.
A key result here is due to \cite{bentkus2003dependence,bentkus2005lyapunov}: Let $\calB$ and $\calC$ be the class of closed balls in $\R^p$ and the class of Borel measurable convex sets in $\R^p$, respectively. Then, when $X,X_1,\dots,X_n$ are i.i.d. with $\E[XX^T]=I_p$ and $\calA = \calB$ or $\calC$, 
\[
\sup_{A \in \calA} | \Pr (S_n \in A) - \Pr (N(0,I_p) \in A)| \lesssim \frac{a(\calA) \E[\| X \|^3]}{\sqrt{n}}
\]
up to a universal constant, where $a(\calB) = 1$ and  $a(\calC) = p^{1/4}$. See \cite{raivc2019multivariate} for an explicit constant in $\lesssim$ when $\calA=\calC$. 
The factor $a(\calA)$ accounts for Gaussian anticoncentration estimates for the class of sets $\calA$; cf. \cite{ball1993reverse} for the class of convex sets. Often, $\E[\| X \|^3]$ scales as $p^{3/2}$ (e.g. consider the case where the coordinates of $X$ are uniformly bounded), so in such cases, for Bentkus' bound to converge to zero, we require $p=o(n^{1/3})$ for $\calA =\calB$ and $p=o(n^{2/7})$ for $\calA = \calC$.
Recent work by \cite{fang2020large} succeeded in relaxing these requirements to $p=o((n/\log^2n)^{1/2})$ for $\calA = \calB$ and $p=o((n/\log^2n)^{2/5})$ for $\calA = \calC$. Still, $p$ needs to be sufficiently small compared with $n$ in both cases. See also \cite{zhilova2020nonclassical, zhilova2020new}.
Another line of research concerns bounding the $2$-Wasserstein distance between the law of $S_n$ and $N(0,\Sigma)$; see \cite{zhai2018high,eldan2020clt,fang2019multivariate,courtade2019existence,bonis2020stein} for related results.
\end{remark}

\subsection{Multiplier and Empirical Bootstraps}
\label{sec: bootstrap}

The high-dimensional CLTs discussed in the preceding section show that the sampling distribution of $S_n$ can be approximated by the Gaussian distribution $N(0,\Sigma)$ uniformly over the rectangles, even when $p \gg n$. However, in practice, the approximating Gaussian distribution $N(0,\Sigma)$ is infeasible since the covariance matrix $\Sigma$ is unknown. We shall consider further estimating the $N(0,\Sigma)$ distribution by using bootstrap methods. In this section, for simplicity, we shall focus our attention on the Gaussian multiplier and empirical bootstraps. We will briefly discuss other bootstrap methods at the end of this section. We first formally define the Gaussian multiplier and empirical bootstraps. Recall $\overline{X}_n = n^{-1}\sum_{i=1}^nX_i$.

\begin{itemize}
\item The \textit{Gaussian multiplier bootstrap}  simulates the data-dependent Gaussian distribution $N(0,\hat{\Sigma})$, where $\hat{\Sigma}$ is the empirical covariance matrix. This can be achieved by simulating the conditional distribution (given the data) of
\[
S_{n}^{B} = \frac{1}{\sqrt{n}} \sum_{i=1}^n \xi_i (X_i - \overline{X}_n), 
\]
where $\xi_1,\dots,\xi_n$ are i.i.d. $N(0,1)$ random variables independent of the data. 
\item The \textit{empirical bootstrap} draws an independent sample $X_1^B,\dots,X_n^B$ from the empirical distribution  $P_n = n^{-1}\sum_{i=1}^n \delta_{X_i}$ and constructs 
\[
S_{n}^{B} = \frac{1}{\sqrt{n}} \sum_{i=1}^n (X_i^B - \overline{X}_n). 
\]
\end{itemize}

\noindent In either case, we approximate $N(0,\Sigma)$ or the sampling distribution of $S_n$ by the conditional distribution of $S_{n}^{B}$. 

Consider first the Gaussian multiplier bootstrap. Observe that, conditionally on the data, we have $S_{n}^{B}  \sim N(0,\hat{\Sigma})$ with $\hat{\Sigma} = n^{-1} \sum_{i=1}^n (X_i-\overline{X}_n)(X_i - \overline{X}_n)^T$.
Thus, for the Gaussian multiplier bootstrap, the problem reduces to comparing two Gaussian distributions with covariance matrices $\hat{\Sigma}$ and $\Sigma$. The following Gaussian comparison inequality, taken from Proposition 2.1 in \cite{chernozhukov2019improved}, implies that the Gaussian multiplier bootstrap is consistent over the rectangles provided that $$\max_{1 \le j ,k \le p}| \hat{\Sigma}_{jk} - \Sigma_{jk} | = o_{P}(1/\log^2 p),$$ which can hold under mild moment conditions. 

\begin{proposition}[Gaussian comparison]\label{coro:g-g-comparison}
Let $\Sigma^\ell, \ell =1,2$ be $p \times p$ covariance matrices such that $\min_{1 \le j \le p} \Sigma_{jj}^1 \ge \underline{\sigma}^2$ for some $\underline{\sigma} > 0$, where $\Sigma_{jk}^\ell$ is the $(j,k)$-component of $\Sigma^{\ell}$. Then, 
\[
\sup_{R \in \calR} \Big| \Pr\big (N(0,\Sigma^1) \in R\big) - \Pr\big (N(0,\Sigma^2) \in R\big)\Big | \leq C\sqrt{\Delta} \log p,
\]
where $C$ is a constant that depends only on $\underline{\sigma}$, and $\Delta = \max_{1\leq j,k\leq p}|\Sigma^1_{jk} - \Sigma^2_{jk}|$. 
\end{proposition}

As for the empirical bootstrap,  in view of the fact that the empirical distribution $P_n$ has mean $\overline{X}_n$ and covariance matrix $\hat{\Sigma}$, we may first approximate the conditional distribution of $S_n^{B}$ by $N(0,\hat{\Sigma})$ over the rectangles by applying the high-dimensional CLT, and then use the Gaussian comparison inequality to further approximate $N(0,\hat{\Sigma})$ by $N(0,\Sigma)$.  The following theorem presents finite sample error bounds for the Gaussian multiplier and empirical bootstraps under a sub-exponential condition.  Recall that $\Pr^B$ denotes the conditional probability given $X_1,\dots,X_n$. 

\begin{theorem}[Bootstrap consistency]
\label{thm: bootstrap}
Under the same assumption as Theorem \ref{thm: HDCLT}, for either the Gaussian multiplier or empirical bootstrap,  we have that, with probability at least $1- C\delta_{1,n}$
\[
\sup_{R \in \calR} \Big|\Pr^B\big (S_n^B \in R \big) - \Pr \big(N(0,\Sigma) \in R\big) \Big| \le  C\delta_{1,n},
\]
where $C$ is a constant that depends only on $\underline{\sigma}$ and $\overline{\sigma}$. 
\end{theorem}

The theorem follows from Lemmas 4.5 and 4.6  in \cite{chernozhukov2019improved}.
Precisely speaking, the actual proof of Theorem \ref{thm: bootstrap} for the empirical bootstrap exploits the fact that the empirical bootstrap matches approximately higher-order moments of $S_n$. See \cite{chernozhukov2019improved} for details. 

As in the high-dimensional CLT, the error bound for the bootstrap depends on $p$ only through $\log p$, which allows $p$ to be much larger than $n$ for the bootstrap to be consistent.
Also, the theorem  imposes no restrictions on the correlation structure between the coordinates of the random vectors, allowing $\Sigma$ to be singular.

Like Theorem \ref{thm:hdclt-poly}, the following theorem, adapted from Theorem 2.6 in \cite{chernozhukov2019improved}, covers polynomial moment conditions.

\begin{theorem}[Bootstrap consistency under polynomial moment condition]
\label{thm:boot-poly}
Under the same assumption as Theorem \ref{thm:hdclt-poly}, for either the Gaussian multiplier or empirical bootstrap,  we have that, with probability at least $1-C(\delta_{1,n} \vee \delta_{2,n}^{[q]})$,  
\begin{align*}
\sup_{R \in \mathcal R} \Big|\Pr^B\big (S_n^B \in R \big) - \Pr (N(0,\Sigma) \in R) \Big| 
 \leq C (\delta_{1,n} \vee \delta_{2,n}^{[q]}),
\end{align*}
where $C$ is a constant that depends only on $q$, $\underline{\sigma}$ and $\overline{\sigma}$. 
\end{theorem}

\begin{figure}[htbp]
\label{fig}
\begin{center}
\includegraphics[width=4.5in]{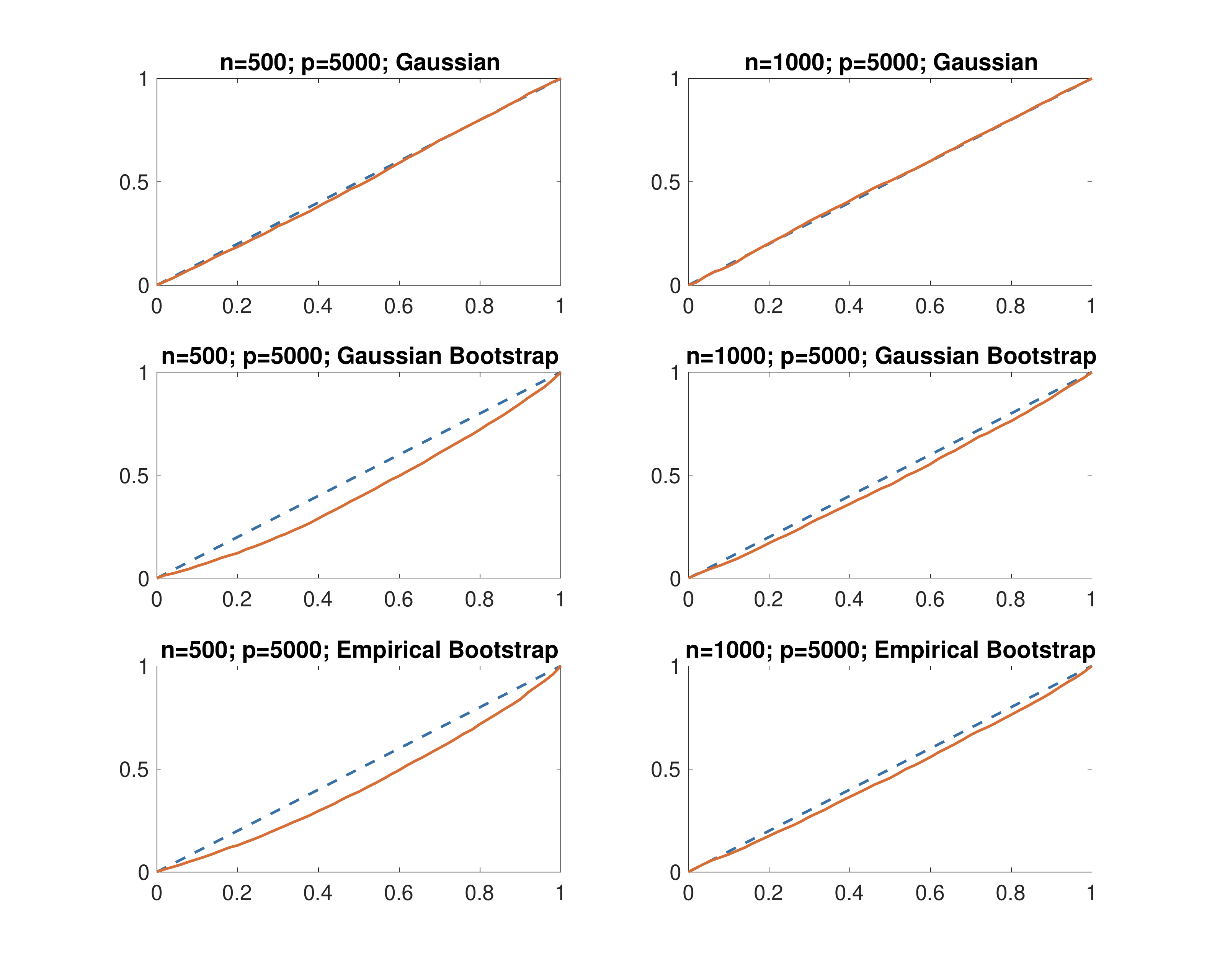}
\caption{P-P plots comparing the distribution of $\| S_n \|_{\infty}$ with its Gaussian, Gaussian multiplier bootstrap, and empirical bootstrap approximations. Here, $X_{i j}$'s are generated as
$X_{i j}=z_{i j}\epsilon_{i}$ with $\epsilon_i$ following the centered (non-symmetric) exponential distribution with parameter $1$, and $z_{i j}$'s are non-stochastic (simulated once using $U[0,1]$ distribution independently across $i$ and $j$).  The dashed line is 45$^\circ$. To generate bootstrap
approximations, we use randomly selected sample of $X_{ij}$'s. The bootstrap distribution provides somewhat conservative approximation in the upper tail in smaller sample, and becomes more accurate in bigger sample. }
\end{center}
\end{figure}

Figure \ref{fig} illustrates the finite sample performance of the Gaussian, Gaussian multiplier bootstrap, and empirical bootstrap approximations for $\| S_n \|_{\infty}$ under a regression setup, namely compares $\Pr(\|S_n\|_\infty \le x)$, $\Pr(\| N(0,\Sigma) \|_{\infty} \le x)$, and $\Pr^B(\| S_n^B \|_\infty \le x)$. The figure indicates
that both Gaussian and bootstrap approximations are reasonably good, especially in the tails.

In applications,  we often normalize the coordinates of the
sample mean by estimates of the standard deviations, so that each coordinate is approximately
distributed as $N(0,1)$. We may estimate the variance $\Sigma_{jj}$ of the $j$-th coordinate of $S_n$ by the sample variance $\hat{\Sigma}_{jj} = n^{-1}\sum_{i=1}^n (X_{ij} - \overline{X}_{n,j})^2$. We shall then consider approximating the sampling distribution of the normalized sample mean $\hat{\Lambda}^{-1/2}S_n$ by the conditional distribution of $\hat{\Lambda}^{-1/2}S_n^B$, where $\hat{\Lambda} = \mathrm{diag} \{ \hat{\Sigma}_{11},\dots,\hat{\Sigma}_{pp} \}$. Let $\Sigma_0 = \Lambda^{-1/2}\Sigma \Lambda^{-1/2}$ denote the correlation matrix of $S_n$. Consider the assumption of Theorem \ref{thm: HDCLT} with $B_n^2 \log^5 (pn) = o(n)$. Then, we have $\max_{1 \le j \le p} |\Sigma_{jj}/\hat{\Sigma}_{jj}-1| = o_{P}(1/\log^2p)$, so that combining the high-dimensional CLT and Gaussian concentration, it holds that $\| (\hat{\Lambda}^{-1/2}-\Lambda^{-1/2})S_n\|_{\infty} = o_{P}(1/\sqrt{\log p})$. Thus, arguing as in the proof of Lemma \ref{lem: ASM}, we have 
\[
 \sup_{R \in \calR} \left|\Pr\big(\hat{\Lambda}^{-1/2}S_n \in R \big) - \Pr \big(N(0,\Sigma_0) \in R\big) \right| \to 0.
\]
Likewise,  for either the Gaussian multiplier or empirical bootstrap, it holds that
\[
 \sup_{R \in \calR} \left|\Pr^B\big(\hat{\Lambda}^{-1/2}S_n^B \in R \big) - \Pr \big(N(0,\Sigma_0) \in R\big) \right| \stackrel{P}{\to} 0.
\]
Similar results hold under polynomial moment conditions. See Appendix A.2 in \cite{chen2019randomized} for relevant arguments. 

The above theoretical results demonstrate that the bootstrap can adequately capture the dependence between the coordinates, thereby yielding asymptotically exact coverage or size control even in high dimensions. Together with the fact that the bootstrap allows for arbitrary correlations between the coordinates, the bootstrap is a particularly powerful inferential tool in high dimensions.

\begin{remark}[Other bootstraps]
The Gaussian multiplier bootstrap is a special case of a wild bootstrap with Gaussian weights. Other bootstrap weights can be also used. 
Indeed, \cite{chernozhukov2019improved} prove Theorem \ref{thm: bootstrap} for more general bootstrap weights, including Mammen's bootstrap \cite{mammen1993bootstrap} as a special case, building on (and improving) important insights in \cite{deng2020beyond} and \cite{koike2021notes}.
\end{remark}

\subsection{Near $n^{-1/2}$ Rates under Structural Assumptions on Covariance Matrix}
\label{sec: parametric rate}

Recall that, in $p=1$, the classical Berry-Esseen theorem shows that, for $X,X_1,\dots,X_n$ i.i.d. univariate $p=1$ random variables with mean zero and unit variance (for simplicity), $\sup_{t \in \R} |\Pr (S_n \le t) - \Pr(N(0,1) \le t)| \lesssim \E[|X|^3]/\sqrt{n}$, and this bound is known to be sharp in terms of dependence on $n$. 
Given this, there has been great interest in deriving high-dimensional CLTs and bootstraps over the rectangles that achieve (near) $n^{-1/2}$ rates while allowing for $p \gg n$. Recent progress shows that such (near) $n^{-1/2}$ rates are possible under structural assumptions on  $\Sigma$. In the following discussion, we assume that $X,X_1,\dots,X_n$ are i.i.d. with mean zero and covariance matrix $\Sigma$. Recall  $S_n = n^{-1/2}\sum_{i=1}^n X_i$. 

\cite{lopes2020bootstrapping} derived a first such result under the variance decay condition. Namely, assuming $\max_{1 \le j \le p}\Var(X_{j}) = O(j^{-a})$ for arbitrarily small $a >0$ and other technical conditions, they derive error bounds of order $n^{-1/2+\delta}$ with arbitrarily small $\delta > 0$ for $\sup_{A \in \calA}|\Pr (S_n \in A) - \Pr (N(0,\Sigma) \in A)|$ for a certain subclass $\calA$ of the rectangles. \cite{lopes2020bootstrapping} also derive similar error bounds for the Gaussian multiplier bootstrap. 

In a different direction, if we assume that the smallest eigenvalue of $\Sigma$ is bounded away from zero, as in \cite{fang2021high}, then it was shown by \cite{chernozhukov2020nearly} that under regularity conditions:
$$\sup_{R \in \calR} |\Pr (S_n \in R) - \Pr(N(0,\Sigma) \in R)| 
\leq C \left (  \frac{B_n^2 (\log p)^3}{n}\right )^{1/2} \log n, $$
and a similar result was obtained for bootstrap approximation. This result builds on and refines a sequence of other important results in this direction obtained by  \cite{fang2021high}, \cite{lopes2020central}, \cite{kuchibhotla2020high}.  Under the nondegeneracy of the covariance matrix $\Sigma$, the above bound implies that the high-dimensional CLT holds if  $\log p = o(n^{1/3})$ up to logarithmic factors, which weakens the previous requirement that $\log p = o(n^{1/5})$. See also \cite{das2021central} who investigate necessary conditions for high-dimensional CLTs.

\section{Selcted Applications}
\label{sec: app}

This section reviews applications of high-dimensional bootstrap to several inference tasks. Specifically, we discuss penalty choice for the Lasso, simultaneous confidence intervals for high-dimensional parameters, estimation and inference for maximum/minumum effects, comparing large covariance matrices, and large-scale multiple testing. There are many other applications of high-dimensional bootstrap, some of which are delineated in the introduction, and the reader is also encouraged  to look at the references therein.

\subsection{Penalty Choice for High-Dimensional Regression}
\hspace*{1mm} Consider a high-dimensional regression with non-Gaussian errors
\[
y_i =x_i^T\beta^{\star} + \epsilon_{i}, \quad \E[\epsilon_i]=0, \quad i=1,\dots,n,
\]
where $y_i$ is a scalar response variable, $x_i$ is a $p$-dimensional vector of (fixed) covariates, and $\epsilon_i$ is an error term.  Here the dimension $p$ of the covariate vector $x_i$ can  be much larger than the sample size $n$, $p \gg n$, but we assume that the model is sparse in the sense that the number of nonzero components of $\beta^{\star}$ is $s = \# \{ j : \beta_j^\star \ne 0 \} \ll n$.  

Arguably, one of the most popular estimates for such high-dimensional linear regression is the Lasso \cite{tibshirani1996regression}, which is defined by
\[
\hat{\beta} \in \argmin_{\beta \in \R^p} \left [ \frac{1}{n} \sum_{i=1}^n (y_i - x_i^{T}\beta )^2 + \lambda \sum_{j=1}^p | \beta_{j} | \right],
\]
where $\lambda \ge 0$ is a tuning parameter. It is well known that the statistical performance of the Lasso crucially relies on the choice of the tuning parameter $\lambda$.
From a seminal work of \cite{bickel2009simultaneous}, for a given confidence level $\alpha \in (0,1)$, if we choose $\lambda$ in such a way that 
\[
\lambda = (1-\alpha)\text{-quantile of } S(\beta^\star),  \quad S(\beta^\star) = 2\cdot \max_{1\leq j \leq p}  \left |    n^{-1} {\textstyle \sum}_{i=1}^n x_{ij} \epsilon_i 
 \right |,
\]
then it holds that $\| \hat{\beta} - \beta^{\star} \|_{2,n} \lesssim \lambda \sqrt{s}$
with probability at least $1-\alpha$, provided that the restricted eigenvalue condition is satisfied; see \cite{bickel2009simultaneous} and also \cite{belloni2013least}. Here $\| \delta \|_{2,n} =\sqrt{ n^{-1}\sum_{i=1}^n (x_i^T \delta)^2 }$.

 If $\epsilon_1,\dots,\epsilon_n$ are i.i.d. sub-Gaussian, then the above choice of $\lambda$ is of order $O(\sqrt{\log p/n})$ (which follows by using a standard concentration inequality) and thus $\| \hat{\beta} - \beta^{\star} \|_{2,n} \lesssim \sqrt{(s\log p)/n}$. However, if the error distribution has heavier tails than sub-Gaussian, then bounding $\lambda$ by concentration inequalities would lead to sub-optimal rates, and the bound itself contains distribution-dependent parameters that are unknown in practice. Instead, we may approximate or estimate $\lambda$ by applying the high-dimensional CLT (homoscedastic case) or using the multiplier bootstrap (heteroscedastic case).

Consider first the case where the error terms are homoscedastic, 
$
\sigma^2 = \E[\epsilon_1^2]=\dots=\E[\epsilon_n^2].
$
Here we assume that $\epsilon_1,\dots,\epsilon_n$ are independent. In this case, we can directly apply the high-dimensional CLT to approximate the quantile of 
$S(\beta^\star)$
 by that of $2\sigma \cdot \max_{1\leq j \leq p}   |    n^{-1}\sum_{i=1}^n x_{ij} \xi_i 
 |$, where $\xi_1,\dots,\xi_n$ are i.i.d. $N(0,1)$. Thus, under homoscedasticity, 
we can approximate $\lambda$ by 
\begin{equation}
\hat{\lambda} = (1-\alpha)\text{-quantile of }\,2\sigma \cdot \max_{1\leq j \leq p}  \left |    n^{-1}{\textstyle \sum}_{i=1}^n  x_{ij} \xi_i \right |.\label{eq: homoscedastic}
\end{equation}
In practice, $\sigma$ is unknown but can be consistently estimated by pre-estimating $\beta^{\star}$ by the Lasso with a crude-choice of the $\lambda$-parameter. 

In general, if the error terms are heteroscedastic, then we first pre-estimate $\beta^{\star}$ (again by using the Lasso with a crude-choice of the $\lambda$-parameter) to construct estimates $\hat{\epsilon}_i$ of $\epsilon_i$, and apply the multiplier bootstrap to estimate the quantile of $\max_{1\leq j \leq p}  |   n^{-1} \sum_{i=1}^n x_{ij} \epsilon_i
|$, namely,
 \begin{equation}
 \hat{\lambda} =\text{conditional $(1-\alpha)$-quantile of $2 \cdot \max_{1\leq j \leq p}  \left |    n^{-1}  {\textstyle \sum}_{i=1}^n x_{ij}\hat{\epsilon}_i \xi_i 
 \right |$}. 
 \label{eq: heteroscedastic}
 \end{equation}
Under regularity conditions, it is shown that the Lasso with this data-driven choice of $\lambda$ satisfies $\| \hat{\beta} - \beta^{\star} \|_{2,n} \lesssim \sqrt{(s \log p)/n}$ with high probability. See Section 4 of \cite{CCK2013AoS} for related results. 
Note that for the justification of both methods (\ref{eq: homoscedastic}) and (\ref{eq: heteroscedastic}), we may apply Theorems \ref{thm:hdclt-poly} and \ref{thm:boot-poly}, which only require that the error terms have finite polynomial moments. 

The function \code{rlasso} in the R package \code{hdm} (\cite{chernozhukov2016high}) implements the above methods of choosing the penalty parameter.\footnote{ The options \code{homoscedastic=TRUE} and \code{X.dependent.lambda=TRUE} in the argument \code{penalty} implement (\ref{eq: homoscedastic}), while \code{homoscedastic=FALSE} (default) and \code{X.dependent.lambda=TRUE} implement (\ref{eq: heteroscedastic}), where the default value of $\alpha$ is set to $\alpha=0.1$} The function \code{rlasso} also offers a joint test of significance of variables -- the test of the null hypothesis that 
$\beta^\star = 0$ -- based on the sup-score statistic $ S(\beta^\star)$ constrained under the null hypothesis $\beta^\star = 0$; the critical value for such test is $\hat \lambda$. 

\subsection{Simultaneous Confidence Intervals for High-Dimensional Parameters}
\label{sec: simultaneous}

Modern statistical and machine learning problems often entail estimation and inference for a large number of parameters, the number of which may exceed the sample size. In such cases, researchers are interested in conducting inference  for not only individual parameters but also groups of parameters simultaneously. Methods of uniform inference for high-dimensional parameters are also a basis of large-scale multiple testing; cf. Section \ref{sec: multiple testing}. 

In this section,  we consider constructing simultaneous confidence intervals (rectangles) for a high-dimensional parameter $\theta^{\star} \in \R^p$. In many settings, e.g., those in examples discussed in below, we have an estimator $\hat{\theta}_n$ for $\theta^{\star}$ such that it admits an asymptotically linear expansion of the form 
\begin{equation}
\hat{\theta}_n - \theta^{\star} = \frac{1}{n}\sum_{i=1}^n \psi_i +\mathsf{R}_n, \label{eq: linear}
\end{equation}
where $\psi_1,\dots,\psi_n$ are independent random vectors (influence functions)  with mean zero and $\mathsf{R}_n$ is a remainder term such that $\| \sqrt{n} \mathsf{R}_n \|_{\infty} = o_{P}(1/\sqrt{\log p})$.
 Applying the high-dimensional CLT to the leading term $n^{-1}\sum_{i=1}^n \psi_i$ yields that (cf. Lemma \ref{lem: ASM}), for $\Sigma = n^{-1} \sum_{i=1}^n \E[\psi_i \psi_i^T]$, 
\[
\sup_{R \in \calR} \Big |\Pr \big(\sqrt{n}(\hat{\theta}_n - \theta^{\star}) \in R\big) - \Pr \big(N(0,\Sigma) \in R\big) \Big | \to 0,
\]
even when $p \gg n$, provided that certain moment conditions on $\psi_1,\dots,\psi_n$ are satisfied. 

We shall estimate the $N(0,\Sigma)$ distribution by using the multiplier bootstrap. In practice, influence functions $\psi_1,\dots,\psi_n$ may be unknown, so we assume that there are suitable estimates $\hat{\psi}_1,\dots,\hat{\psi}_n$ for the influence functions. We then apply the multiplier bootstrap to the estimated influence functions, i.e.,
\[
S_n^B = \frac{1}{\sqrt{n}} \sum_{i=1}^n \xi_i (\hat{\psi}_i - \overline{\hat{\psi}}),
\]
where $ \overline{\hat{\psi}} = n^{-1}\sum_{i=1}^n \hat{\psi}_i$. Let $\hat{\Sigma} =  n^{-1}\sum_{i=1}^n (\hat{\psi}_{i} -  \overline{\hat{\psi}})(\hat{\psi}_{i} -  \overline{\hat{\psi}})^T$.  In view of the Gaussian comparison inequality (Proposition \ref{coro:g-g-comparison}), the Gaussian multiplier bootstrap is consistent over the rectangles,
$
\sup_{R \in \calR}  | \Pr^B (S_n^B \in R) - \Pr (N(0,\Sigma) \in R) | \stackrel{P}{\to} 0,
$
provided that $\max_{1 \le j,k \le p} | \hat{\Sigma}_{jk} - \Sigma_{jk}| = o_{P}(1/\log^2 p)$,
 which can hold even when $p \gg n$. Alternatively, we may apply the empirical bootstrap by generating an independent sample $\hat{\psi}_1^B,\dots,\hat{\psi}_n^B$ from the empirical distribution $n^{-1}\sum_{i=1}^n \delta_{\hat{\psi}_i}$ and construct $S_n^B = n^{-1/2}\sum_{i=1}^n (\hat{\psi}_i^B -  \overline{\hat{\psi}})$. 
 
 Also, the high-dimensional CLT and bootstrap consistency hold for the normalized statistics $\hat{\Lambda}^{-1/2}\sqrt{n}(\hat{\theta}_n - \theta^{\star})$ and $\hat{\Lambda}^{-1/2}S_n^B$, under regularity conditions, where  $\hat{\Lambda} = \mathrm{diag} \{ \hat{\Sigma}_{11},\dots,\hat{\Sigma}_{pp} \}$ with $\hat{\Sigma}_{jj} = n^{-1}\sum_{i=1}^n (\hat{\psi}_{ij} - \overline{\hat{\psi}}_{j})^2$ (see the discussion at the end of Section \ref{sec: bootstrap}), namely,
 \begin{equation}
 \label{eq: consistency}
 \begin{split}
&\sup_{R \in \calR} \Big|\Pr \big(\hat{\Lambda}^{-1/2}\sqrt{n}(\hat{\theta}_n - \theta^{\star}) \in R\big) - \Pr \big(N(0,\Sigma_0) \in R\big) \Big| \to 0 \quad \text{and}\\
&\sup_{R \in \calR} \Big|\Pr^B \big(\hat{\Lambda}^{-1/2}S_n^B \in R\big) - \Pr \big(N(0,\Sigma_0\big) \in R)\Big | \stackrel{P}{\to} 0 . 
\end{split}
 \end{equation}
 Here $\Sigma_0 = \Lambda^{-1/2} \Sigma \Lambda^{-1/2}$ with $\Lambda = \mathrm{diag} \{ \Sigma_{11},\dots,\Sigma_{pp} \}$, and $p$ is allowed to increase with $n$, $p=p_n \to \infty$ (we assume that the diagonal elements $\Sigma_{jj}$ are bounded away from zero). 
For a given $\alpha \in (0,1)$, let 
\[
\hat{q}_{1-\alpha} = \text{conditional $(1-\alpha)$-quantile of $\| \hat{\Lambda}^{-1/2}S_n^{B} \|_{\infty}$}.
\]
We claim that the rectangle of the form
\[
\hat{R}_{1-\alpha}= \prod_{j=1}^p \left [\hat{\theta}_{n,j} \pm \hat{\Sigma}_{jj}^{1/2}\hat{q}_{1-\alpha}/\sqrt{n}  \right ] =: \prod_{j=1}^p \hat{\mathrm{CI}}_{j,1-\alpha}
\]
contains $\theta^{\star}$ with probability $1-\alpha + o(1)$. Indeed, $\theta^{\star} \in \hat{R}_{1-\alpha}$ if and only if $\| \hat{\Lambda}^{-1/2}\sqrt{n}(\hat{\theta}_n - \theta^{\star}) \|_{\infty} \le \hat{q}_{1-\alpha}$. Let $q_{1-\alpha}$ denote the $(1-\alpha)$-quantile of $\| N(0,\Sigma_0) \|_{\infty}$. Then, from the bootstrap consistency result (\ref{eq: consistency}), there exists a sequence $\epsilon_n \to 0$ such that $q_{1-\alpha-\epsilon_n} \le \hat{q}_{1-\alpha} \le q_{1-\alpha+\epsilon_n}$ with probability approaching one (cf. Theorem 2.5 in \cite{belloni2018high}). Combining the high-dimensional CLT from (\ref{eq: consistency}) and the fact that $\|  N(0,\Sigma_0) \|_{\infty}$ has a continuous distribution function (cf. Proposition \ref{lem: Nazarov}), we see that
\[
\Pr (\| \hat{\Lambda}^{-1/2}\sqrt{n}(\hat{\theta}_n - \theta^{\star}) \|_{\infty} \le \hat{q}_{1-\alpha})  \le 1-\alpha+ \epsilon_n + o(1) = 1-\alpha + o(1). 
\]
The reverse inequality follows similarly. Conclude that $\hat{R}_{1-\alpha}$ is a valid simultaneous confidence interval for $\theta^{\star}$ with level $1-\alpha+o(1)$, i.e., $\Pr (\theta^{\star} \in \hat{R}_{1-\alpha}) = 1-\alpha+o(1)$. 

\begin{example}[Randomized control trials with many outcomes]
Consider a randomized control trial with $n$ participants, where each participant $i=1,\dots,n$ is randomly assigned to either the treatment group ($D_i = 1$) or the control group ($D_i = 0$).
For each participant $i$, we observe a large number of outcome variables $Y_{i} = (Y_{i1},\dots,Y_{ip})^T \in \R^p$. 
Let $\gamma$ denote the probability of being assigned to the treatment group. Then, the average treatment effect for outcome variable $Y_{ij}$ is defined as $\theta_j^{\star} = \E[Y_{ij} \mid D_i=1] - \E[Y_{ij} \mid D_i=0]$,
which can be estimated as
\[
\hat{\theta}_{n,j} = \frac{1}{n}\sum_{i=1}^n \left(\frac{D_iY_{ij}}{\gamma} - \frac{(1-D_i)Y_{ij}}{1-\gamma} \right ).
\]
In this case, the estimator is exactly linear (i.e., $\mathsf{R}_n = 0$ in  (\ref{eq: linear})), and the method above is readily applicable to conducting simultaneous inference for $\theta^{\star} = (\theta_{1}^{\star},\dots,\theta_{p}^{\star})^T$. 
\end{example}

\begin{example}[Lasso-penalized  $M$-estimation]
For $M$-estimation problems involving high-dimensional parameters, such as high-dimensional (generalized) linear models or quantile regression models, one can construct estimators admitting asymptotically linear forms (like (\ref{eq: linear})) by ``debiasing'' regularized estimators or constructing orthogonal estimating equations \cite{zhang2014confidence,belloni2014inference,van2014asymptotically,javanmard2014confidence,Belloni2015,belloni2018uniformly}. \cite{Belloni2015}  apply the multiplier bootstrap to debiased lasso estimators to  construct simultaneous confidence intervals for the slope vector in a high-dimensional median regression and other $M$-regression problems, viewed through the lenses of $Z$-estimation framework for many target parameters. They develop orthogonal estimating equations that give rise to estimators permitting asymptotically linear forms of the type given in equation (\ref{eq: linear}), for which they apply the multiplier bootstrap to construct simultaneous confidence intervals. 
\cite{ZhangCheng2017JASA} and \cite{Dezeure2017} apply the multiplier bootstrap to the debiased Lasso to  construct simultaneous confidence intervals for the slope vector in a high-dimensional linear mean regression model.
See also \cite{belloni2018uniformly} for the extension to models involving $p$ functional parameters, with $p \gg n$.

Finally, the function \code{rlassoEffects} in the R package \code{hdm} implements a form of debiased Lasso based on \cite{belloni2014inference} and the multiplier bootstrap based on \cite{Belloni2015} to construct simultaneous confidence intervals for linear mean regression, where \code{confint} extracts results on confidence intervals.  

\end{example}

\begin{example}[Post-selection inference]
Simultaneous confidence intervals can be adapted to post model selection inference. Suppose that $\{ 1,\dots, p \}$ now represents the set of models and $\theta_j^{\star}$ corresponds to a parameter for model $j$. Researchers select a model $\hat{j}$ using the data and construct a confidence interval for $\theta_{\hat{j}}^{\star}$. However, a naive plug-in confidence interval, i.e., $\big[ \hat{\theta}_{n,\hat{j}} \pm \frac{\hat{\Sigma}_{\hat{j},\hat{j}}z_{1-\alpha/2}}{\sqrt{n}} \big]$ (with $z_{1-\alpha/2}$ being the $(1-\alpha/2)$-quantile of $N(0,1)$), will not be a valid confidence interval for $\theta_{\hat{j}}^{\star}$ since this construction ignores the randomness involved in the model selection step $\hat{j}$. 

To tackle this problem, \cite{kuchibhotla2020valid} consider replacing the normal critical value $z_{1-\alpha/2}$ with $\hat{q}_{1-\alpha}$ above (see also \cite{berk2013valid}), and show that the resulting confidence interval, $\hat{\mathrm{CI}}_{\hat{j},1-\alpha} =\big [ \hat{\theta}_{n,\hat{j}} \pm \frac{\hat{\Sigma}_{\hat{j},\hat{j}}\hat{q}_{1-\alpha}}{\sqrt{n}} \big]$, satisfies the desired coverage requirement. Indeed, by construction, we see that
\[
\Pr \big(\theta_{\hat{j}}^{\star} \in \hat{\mathrm{CI}}_{\hat{j},1-\alpha} \big) \ge \Pr \Big (\bigcap_{j=1}^p \{ \theta_j^{\star} \in \hat{\mathrm{CI}}_{j,1-\alpha} \} \Big) 
= \Pr \big(\theta^{\star} \in \hat{R}_{1-\alpha} \big) = 1-\alpha+o(1). 
\]
The preceding argument immediately extends to the case where each $\theta_j^{\star}$ is a vector of increasing dimensions (in that case, each $\hat{\mathrm{CI}}_{j,1-\alpha}$ will be a rectangle). See \cite{kuchibhotla2021post} for a review on post-selection inference. 
\end{example}

\subsection{Estimation and Inference for Maximum/Minimum Effects}

In certain applications, researchers are interested in the maximum or minimum of high-dimensional parameters.  We start with discussing two such examples. 

\begin{example}[Intersection bounds]
Suppose that the parameter of interest $\vartheta^{\star}$ is known to lie  within the intervals $[\theta^l_j,\theta^u_j]$ for $j=1,\dots,p$, i.e.,
\[
\max_{1 \le j \le p} \theta_j^l \le \vartheta^{\star} \le \min_{1 \le j \le p}\theta_j^u,
\]
where lower and upper bounding parameters $\theta_j^l, \theta_j^u$ for $j=1,\dots,p$ are estimable from data, but no other information is available to estimate $\vartheta^{\star}$. That is, the parameter $\vartheta^\star$ is identified only up to the set $\big [ \max_{1 \le j \le p} \theta_j^l, \min_{1 \le j \le p}\theta_j^u\big ]$. Such partially identified models commonly appear in econometrics, see for example, \cite{manski2010partial}, \cite{manski2009more} and \cite{chesher2017generalized}, where instrumental variables are used to bound average causal effects in the presence of missing data or latent confounders. \cite{chernozhukov2013intersection} develop a generic method of inference for the bounding parameters and establish the asymptotic validity allowing $p$ to increase with $n$. 
\end{example}

\begin{example}[Best subgroup effect]
Statistical analysis of treatment heterogeneity across different subgroups has received increasing attention. 
In clinical trials, a new treatment may be only marginally effective for the overall population, but often appears to be promising to certain subgroups. 
Motivated by this, \cite{guo2021inference} consider inference for the best subgroup effect $\theta^{\star}_{\max} =\max_{1 \le j \le p} \theta_j^{\star}$ among $p$ (possibly overlapping) groups, where $\theta_j^{\star}$ is the treatment effect for the $j$-th group for $j=1,\dots,p$. \cite{guo2021inference} propose a modified bootstrap confidence bound for $\theta^{\star}_{\max}$ and a method of bias correction for the plug-in estimator, and establish their asymptotic validity under the setting that $p$ (the number of subgroups) is fixed. 
\end{example}

Consider the setting of Section \ref{sec: simultaneous} and estimation of $\vartheta^{\star} = \max_{1 \le j \le p}\theta_j^\star$ (the minimum effect can be dealt with analogously, as $\min_{1 \le j \le p}\theta_j^\star = - \max_{1 \le j \le p} (-\theta_j^\star)$). As observed in \cite{chernozhukov2013intersection} and \cite{guo2021inference}, the plug-in estimator $\max_{1 \le j \le p}\hat{\theta}_{n,j}$ tends to be upward biased in the finite sample. To tackle this issue, \cite{chernozhukov2013intersection} proposed the following precision corrected estimator 
\[
\hat{\vartheta}_n (\alpha) = \max_{1 \le j \le p} \big [ \hat{\theta}_{n,j} - \hat k_{1-\alpha} \hat{\Sigma}^{1/2}_{jj}/\sqrt{n} \big ],
\]
which can make the estimator upward $(1-\alpha)$-quantile unbiased, for example median upward unbiased for $\alpha=1/2$, and  $\hat k_{1-\alpha}$
is the bootstrap estimate of  $(1-\alpha)$-quantile of 
the studentized maximum estimation error $ T= \sqrt{n} \max_{1 \le j \le p}  (\hat \theta_{n,j} - \theta_j^{\star})/\hat{\Sigma}^{1/2}_{jj}.$ Indeed,
\[
\begin{split}
\Pr \big ( \hat{\vartheta}_n (\alpha) \le \vartheta^\star \big ) &\ge \Pr \Big ( \hat{\theta}_{n,j} - \theta_j^\star \le \hat{k}_{1-\alpha} \hat{\Sigma}^{1/2}_{jj}/\sqrt{n} \ \text{for all $j=1,\dots,p$} \Big) \\
&= \Pr (T \le \hat{k}_{1-\alpha}) = 1-\alpha + o(1), 
\end{split}
\]
under regularity conditions like those in Section \ref{sec: simultaneous}. Thus, the precision correction guarantees that the estimator is upward biased with probability at most $\alpha+o(1)$. Also, $[\hat{\vartheta}_n(\alpha),\infty)$ is an asymptotically valid one-sided confidence interval for $\vartheta^*$ with level $1-\alpha+o(1)$. 
The approach can be refined by a conservative pre-estimation of the argmax set $J_0 = \arg \max_{1 \leq j \leq p}\theta^*_j$, and then working with such set in place of $\{1,\dots, p\}$.

\begin{example}[Best policy estimation]
The problem of estimating the argmax $J_0$ set brought above is interesting in its own right, and arises in optimal policy analysis, where $j$'s correspond to different policies and $\theta^\star_j$'s correspond  to gains (average treatment effects) associated to those policies over the status quo (control) policy; see, e.g. \cite{athey:wager}. In this case $J_0$ is a set of best policies. \cite{chernozhukov2013intersection} propose to estimate $J_0$ by 
$$
\hat J =\Big \{ j \in \{1,\dots,p \}: \hat \theta_{n,j} + \hat  q_{1-\beta} \hat{\Sigma}^{1/2}_{jj}/\sqrt{n} 
\geq \max_{1 \leq j \leq p}  [\hat \theta_{n,j} - \hat q_{1-\beta}\hat{\Sigma}^{1/2}_{jj}/\sqrt{n}] \Big \},
$$
where $\hat q_{1-\beta}$ is the bootstrap estimate of  $(1-\beta)$-quantile of  the studentized maximum absolute estimation error $ \bar T= \sqrt{n} \max_{1 \le j \le p} | \hat \theta_{n,j} - \theta_j^{\star}|/\hat{\Sigma}^{1/2}_{jj}.$ Then $ J_0 \subset \hat J $ with probability no less than $1- \beta - o(1)$ under regularity conditions like those in Section \ref{sec: simultaneous}.

\end{example}

\subsection{Comparing Large Covariance Matrices}
In modern genomics, understanding how the dependencies among many genes vary between different biological states (e.g., healthy or with disease) has received significant attention.  Formally, the problem amounts to comparing large covariance matrices across different populations.  Let $X$ and $Y$ be two random vectors in $\R^p$ with covariance matrices 
$$
\Sigma_1 = \Big (\sigma_{1,jk} \Big )_{1 \le j,k \le p} \  \ \mathrm{ and  }  \ \ \Sigma_2 =  \Big (\sigma_{2,jk} \Big )_{1 \le j,k \le p}, 
$$
respectively,
and  consider testing the hypothesis $H_0: \Sigma_1 = \Sigma_2$ against the alternative $H_1: \Sigma_1 \ne \Sigma_2$ in high dimensions. 

Let $X_1,\dots,X_n$ and $Y_1,\dots,Y_m$ be independent observations from $X$ and $Y$, respectively, and let $\hat{\Sigma}_1 = (\hat{\sigma}_{1,jk})_{1 \le j,k \le p}$ and $\hat{\Sigma}_2 = (\hat{\sigma}_{2,jk})_{1 \le j,k \le p}$ be the corresponding empirical covariance matrices, respectively (e.g., $\hat{\Sigma}_1 = n^{-1}\sum_{i=1}^n (X_i - \overline{X}_n)(X_i - \overline{X}_n)^T$ with $\overline{X}_n = n^{-1}\sum_{i=1}^n X_i$). \cite{chang2017} consider the max-type test statistic $\hat{T}_{\max} = \max_{1 \le j \le k \le p} |\hat{t}_{jk}|$,
where 
\[
\hat{t}_{jk} = \frac{\hat{\sigma}_{1,jk}-\hat{\sigma}_{2,jk}}{\sqrt{n^{-1}\hat{s}_{1,jk}+m^{-1}\hat{s}_{2,jk}}}.
\]
Here $\hat{s}_{1,jk} = n^{-1}\sum_{i=1}^n \{ (X_{ij} -\overline{X}_{n,j})(X_{ik} -\overline{X}_{n,k}) - \hat{\sigma}_{1,jk} \}^2$ and $\hat{s}_{2,jk}$ is defined analogously. 

To calibrate critical values for the test, \cite{chang2017} propose the following Gaussian multiplier bootstrap procedure:

\begin{enumerate}
\item Generate i.i.d. $N(0,1)$ random variables $\xi_1,\dots,\xi_{n+m}$ independent of the data and construct the multiplier bootstrap statistic $\hat{T}_{\max}^\dagger = \max_{1 \le j \le k \le p} |\hat{t}_{jk}^{\dagger}|$,
where 
\[
\hat{t}_{jk}^\dagger = \frac{\hat{\sigma}_{1,jk}^\dagger - \hat{\sigma}_{2,jk}^\dagger }{\sqrt{n^{-1}\hat{s}_{1,jk}+m^{-1}\hat{s}_{2,jk}}}
\]
with $\hat{\sigma}_{1,jk}^\dagger = n^{-1}\sum_{i=1}^n \xi_i \{(X_{ij} -\overline{X}_{n,j})(X_{ik} -\overline{X}_{n,k}) - \hat{\sigma}_{1,jk}\}$ and $\hat{\sigma}_{2,jk}^\dagger$ defined analogously. 
\item For a given $\alpha \in (0,1)$, compute the critical value $\hat{q}_{1-\alpha}$ as the conditional $(1-\alpha)$-quantile of $\hat{T}_{\max}^\dagger$. 
\end{enumerate}
Assuming that $m$ is of comparative size as $n$, \cite{chang2017} show that the test has an asymptotic size $\alpha$ even if $p \gg n$, namely, they show that
$\Pr (\hat{T}_{\max}^\dagger > \hat{q}_{1-\alpha}) = \alpha+o(1)$
if $H_0$ holds true.  (From the results reviewed in this article, the improved condition on $p$ is that $\log p = o(n^{1/5})$.) Further, \cite{chang2017} combine this test with the Benjamini–Hochberg procedure to develop a gene clustering algorithm.

\subsection{Large-Scale Multiple Testing}
\label{sec: multiple testing}

Suppose that we have a large collection of null $H_j$ and alternative $H_j'$ hypotheses for $j=1,\dots,p$. Such large-scale multiple testing problems commonly appear in biological applications. 
Here, we are interested in testing these hypotheses simultaneously for all $j=1,\dots,p$, so we aim  to construct a testing procedure that would reject at least one
true null hypothesis with probability at most $\alpha+o(1)$, uniformly
over the set of true null hypotheses.  Procedures with this property are said to control the Family-Wise Error Rate (FWER).

To this end, one can adopt the step-down procedure of \cite{romano2005}. 
Specifically, consider the setting of Section \ref{sec: simultaneous} and the testing problems $H_j: \theta_j^\star \le 0$ against $H_j': \theta_j^\star > 0$ for $j=1,\dots,p$ (the equality case $\theta_j^\star = 0$ can be dealt with by splitting $\theta_j^\star = 0$ into two hypotheses $\theta_j^\star \le 0$ and $-\theta_j^\star \le 0$). For each $j=1,\dots,p$, let $\hat{t}_j = \sqrt{n}\hat{\theta}_{n,j}/\hat{\Sigma}_{jj}^{1/2}$ be the studentized test statistic for $H_j$ against $H_j'$. 

 For a subset $w\subset \{1,\dots,p\}$, let $c_{1-\alpha,w}$ be the bootstrap estimate of
the $(1-\alpha)$-quantile of $\max_{j\in w}\sqrt{n}(\hat{\theta}_{n,j}-\theta_j^\star)/\hat{\Sigma}_{jj}^{1/2}$ by either applying the Gaussian multiplier or empirical bootstrap (cf. Section \ref{sec: simultaneous}). On the first step, let $w(1)=\{1,\dots,p\}$. Reject
all hypotheses $H_{j}$ satisfying $\hat{t}_{j}>c_{1-\alpha,w(1)}$.
If no null hypothesis is rejected, then stop. If some $H_{j}$ are
rejected, let $w(2)$ be the set of all null hypotheses that
were not rejected in the first step. In step $l\geq2$, let $w(l)\subset \{1,\dots,p \}$
be the subset of null hypotheses that were not rejected up to step
$l$. Reject all hypotheses $H_{j}$, $j\in w(l)$, satisfying $\hat{t}_{j}>c_{1-\alpha,w(l)}$.
If no null hypothesis is rejected, then stop. If some $H_{j}$ are
rejected, let $w(l+1)$ be the subset of all null hypotheses
among $j\in w(l)$ that were not rejected. Proceed in this way until
the algorithm stops. 

\cite{CCK2013AoS} show that this step-down procedure can achieve the FWER control under regularity conditions, even when $p \gg n$, extending the analysis of \cite{romano2005} to high dimensions.  See Section 5 in \cite{CCK2013AoS} and also Section 2.4 in \cite{belloni2018high} for more details.

The following R code illustrates implementation of Romano--Wolf's step-down procedure for multiple one-sample $t$-tests. We compute the adjusted $p$-values following \cite{RoWo16}. The dataset \texttt{Fund} is taken from the \texttt{ISRL2} package associated with the textbook \cite{JWHT21}. The dataset contains $n=50$ rows and $p=2000$ columns, and each column corresponds to the returns of a hedge fund manager. Writing $\mu_j$ for the $j$-th fund manager's mean return, we test for $H_j:\mu_j\leq0$ against $H_j':\mu_j>0$ for all $j=1,\dots,2000$. See Section 13.3 of \cite{JWHT21} for more illustration of this dataset. 

\begin{lstlisting}
library(boot); library(ISLR2)
n <- nrow(Fund); y <- scale(Fund) 
m <- attr(y, "scaled:center"); s <- attr(y, "scaled:scale")
tstat <- sqrt(n) * m/s; ord <- order(tstat) 
# Empirical bootstrap
R <- 499 # Number of bootstrap replicates
mystat.e <- function(y, i){
  tstar <- sqrt(n) * colMeans(y[i, ]) # bootstrapped t statistics
  tmax <- cummax(tstar[ord]); return(tmax >= tstat[ord])
}
set.seed(111)
res.e <- boot(y, mystat.e, R) # bootstrap
count.e <- colSums(res.e$t) 
pval.e <- double(ncol(Fund))
pval.e[rev(ord)] <- (1 + cummax(rev(count.e)))/(R + 1)
which(pval.e < 0.1) # Two funds are detected
# Multiplier bootstrap
mystat.m <- function(y){
  tstar <- sqrt(n) * colMeans(y) # bootstrapped t statistics
  tmax <- cummax(tstar[ord]); return(tmax >= tstat[ord])
}
ran.gen <- function(ystar, mle) rnorm(length(ystar)) * ystar
res.m <- boot(y, mystat.m, R, sim = "parametric", ran.gen = ran.gen)
count.m <- colSums(res.m$t)
pval.m <- double(ncol(Fund))
pval.m[rev(ord)] <- (1 + cummax(rev(count.m)))/(R + 1)
which(pval.m < 0.1) # Two funds are detected
\end{lstlisting}

Finally, the R package \code{hdm} offers functionality on multiple hypothesis testing in high-dimensional approximately sparse linear regression models based upon Romano--Wolf step down procedures (see \cite{bach2018valid} for documentation).

\section{Concluding Remarks}
\label{sec: discussion}

The field of high-dimensional bootstrap has seen a rapid development, and this article reviewed the main ideas and key techniques used there. 
Notably, the bootstrap offers the following advantages to inference for high-dimensional data:

\begin{itemize}
    \item The bootstrap consistency holds in high dimensions ($p \gg n$) without relying on explicit limit distributions. In some situations, one can find limit distributions of one-dimensional functionals of high-dimensional data, but the derivation typically requires restrictive conditions such as weak dependence across the coordinates. The justification of the bootstrap does not require such restrictive assumptions. 
    \item The bootstrap is able to automatically capture the dependence structure among the coordinates, thereby yielding asymptotically exact inference methods. This is so even when the data dimension is much larger than the sample size, where the data coordinates possess a potentially complex dependence structure. 
\end{itemize}

\noindent The theoretical development stimulated many new applications of bootstrap to high-dimensional inference tasks, some of which we reviewed in this article. 

We end this article with commenting on a couple of future research topics. First, Section \ref{sec: bootstrap} presents error bounds for the Gaussian multiplier and empirical bootstraps, and the same error bounds hold for other wild bootstraps such as Mammen's bootstrap. As observed numerically in \cite{deng2020beyond} and \cite{chernozhukov2019improved}, however, the empirical and Mammen's bootstraps perform (slightly) better than the Gaussian multiplier bootstrap in practice. Deriving error bounds that certify a specific bootstrap to be preferred over others is an interesting direction in future research. Also, while there is considerable progress on extending high-dimensional bootstrap for (temporally or graph-) dependent data, the derived error bounds in the dependent case are substantially slower than the independent case. More research is needed to explore sharp error bounds for high-dimensional bootstrap for dependent data. Finally, the literature on high-dimensional bootstrap has been focused on the case where the approximating distribution is Gaussian. However, several important statistics such as degenerate $U$-statistics have non-Gaussian approximating distributions. With an exception of \cite{koike2019mixed}, extending the scope of high-dimensional bootstrap to non-Gaussian approximating distributions  is open.

\bibliographystyle{amsalpha}

\bibliography{biblio}

\end{document}